\documentclass[leqno,11pt]{article}%
\usepackage{amsmath,amssymb}
\usepackage{amsthm}
\usepackage{dsfont}
\usepackage{amsfonts}%
\usepackage{color}
\usepackage[utf8]{inputenc}

\newtheorem{Theorem}{\hspace{\parindent}\bf Theorem}[section]
\newtheorem{Lemma}{\hspace{\parindent}\bf Lemma}[section]
\newtheorem{Proposition}{\hspace{\parindent}\bf Proposition}[section]

\newtheorem{Remark}{\hspace{\parindent}\bf Remark}[section]

\begin{document}

\title{\textbf{Coupling local and nonlocal evolution equations}}
\author{by\\Alejandro G\'{a}rriz, Fernando Quir\'os, and Julio D. Rossi}

\maketitle

\begin{abstract}
We prove existence, uniqueness and several qualitative properties for evolution equations that combine
local and nonlocal diffusion operators acting in different subdomains and coupled in such a way that
the resulting evolution equation is the gradient flow of an energy functional. We deal with the Cauchy, Neumann and Dirichlet problems, in the last two cases with zero boundary data. For the first two problems we prove that the model preserves the total mass. We also study the behaviour  of the solutions for large times. Finally, we show that we can recover the usual heat equation (local diffusion) in a limit procedure when we rescale the nonlocal kernel in a suitable way.
\end{abstract}


\noindent{\makebox[1in]\hrulefill}\newline2010 \textit{Mathematics Subject
Classification.} 
35K55, 
35B40, 
35A05. 

\section{Introduction and main results}
\label{sect-introduction} \setcounter{equation}{0}

If you think about a linear diffusion equation, probably the first one that will come to your  mind is the classical heat equation
\begin{equation} \label{heat.equation}
u_t = \Delta u.
\end{equation}
This equation is naturally associated with the energy
\begin{equation}
\label{eq:energy.heat}
E(u) = \int \frac{|\nabla u|^2}{2},
\end{equation}
in the sense that \eqref{heat.equation} is the gradient flow associated to $E(u)$, see \cite{Evans}.

If you go one step further and consider nonlocal diffusion problems,
one popular choice is
\begin{equation} \label{eq:nonlocal.equation}
u_t (x,t) = \int_{\mathbb{R}^N} J(x-y)(u(y,t)-u(x,t))
\, {\rm d}y,
\end{equation}
where $J: \mathbb{R}^N \to \mathbb{R}$ is a nonnegative, radial function with
$\int_{\mathbb{R}^N} J =1$. Notice that the diffusion of the density $u$ at a point $x$ and time $t$ depends on the values of $u$ at all points in the set $x+\mathop{\rm supp\, }J$, which is what makes the diffusion operator nonlocal. Evolution equations of this form
and variations of it have been recently widely used to model
diffusion processes, see \cite{BCh, BFRW, CF, C, CERW, delia1, F, FW, S, W, Z}. As stated in \cite{F}, if $u(x,t)$ is thought of
as the density of a single population at the point $x$ at time $t$, and
$J(x-y)$ is regarded as the probability distribution of jumping from location
$y$ to location~$x$, then  the
rate at which individuals are arriving to position $x$ from all other places is given by $\int_{\mathbb{R}^N} J(y-x)u(y,t)\,{\rm d}y$, while  the rate at which they are
leaving location $x$ to travel to all other sites is given by $-\int_{\mathbb{R}^N} J(y-x)u(x,t)\, {\rm d}y=-u(x,t)$. Therefore,  in the
absence of external or internal sources, the
density $u$ satisfies equation (\ref{eq:nonlocal.equation}).
In this case there is also an energy that governs the
evolution problem, namely
\begin{equation}
\label{eq:energy.nonlocal}
E(u) = \frac14 \iint J(x-y) (u(y) - u(x))^2 \, {\rm d}x{\rm d}y.
\end{equation}

In the present paper we consider an energy which is local in certain subdomain and nonlocal in the complement, and study the associated gradient flow.
We will show that the Cauchy, Neumann and Dirichlet problems, in the last two cases with zero boundary data, for this equation are well posed. Moreover, we will prove that the solutions to these problems share several properties with their local and nonlocal counterparts~\eqref{heat.equation} and~\eqref{eq:nonlocal.equation}: conservation of mass,  for the Cauchy and Neumann problems, comparison principles, and asymptotic behaviour as $t\to \infty$.

\subsection{The Cauchy problem}
Let $\Gamma$ be a smooth hypersurface  that divides the space $\mathbb{R}^N$ in two smooth domains $\Omega_\ell$ and $\Omega_{n\ell}$.
We introduce the energy \begin{equation}\label{def:funcional_cauchy.intro}
E(u):=\underbrace{\int_{\Omega_{\ell}} \frac{|\nabla u|^2}{2}}_{\mathcal{L}(u)} + \underbrace{\frac{\alpha}{2}\int_{\Omega_{n\ell}}\int_{\mathbb{R}^N} J(x-y)(u(y)-u(x))^2\, {\rm d}y{\rm d}x}_{\mathcal{N}(u)},
\end{equation}
where $\alpha$ is a fixed positive constant. Thus, the energy functional has two parts, a local one $\mathcal{L}(u)$,  that resembles the energy functional~\eqref{eq:energy.heat} for the equation~\eqref{heat.equation}, and a nonlocal part, $\mathcal{N}(u)$, similar to the energy~\eqref{eq:energy.nonlocal} associated with the nonlocal heat equation~\eqref{eq:nonlocal.equation}.

We would like our equation to be the the gradient flow of the energy functional~\eqref{def:funcional_cauchy.intro}. To be more precise, $u$ will be the solution of the ODE (in an infinite dimensional space)
 $u'(t)= -\partial E[u(t)]$, $t\ge0$, $u(0)=u_0$, where $\partial E[u]$ denotes the subdifferential of $E$ at the point $u$. To compute the subdifferential, we obtain the derivative of $E$ at $u$ in the direction of $\varphi\in C^\infty_0(\mathbb{R}^N)$,
$$
\begin{aligned}
\displaystyle \partial_\varphi E(u)=& \displaystyle\lim\limits_{h\downarrow 0}\frac{E(u+h\varphi) - E(u)}{h}
\\[6pt]
=& \displaystyle \int_{\Omega_\ell} \nabla u\cdot\nabla \varphi+\alpha\int_{\Omega_{n\ell}}\int_{\mathbb{R}^N} J(x-y)(u(y)-u(x))(\varphi(y)-\varphi(x))\, {\rm d}y{\rm d}x.
\end{aligned}
$$
Thus, if $u$ were smooth, we would have
$$
\begin{aligned}
	\displaystyle \partial_\varphi& E(u)= \int_{\Gamma}\varphi\partial_\eta u -\int_{\Omega_\ell} \left\lbrace \Delta u (x) + \alpha  \int_{\Omega_{n\ell}} J(x-y)(u(y)-u(x))\,{\rm d}y \right\rbrace\varphi(x)\, {\rm d}x
\\[6pt]
 \displaystyle&- \alpha\int_{\Omega_{n\ell}} \varphi(x) \int_{\Omega_\ell} J(x-y)(u(y)-u(x))\, {\rm d}y {\rm d}x
- 2  \alpha\int_{\Omega_{n\ell}}\varphi(x) \int_{\Omega_{n\ell}} J(x-y)(u(y)-u(x))\,{\rm d}y {\rm d}x,
\end{aligned}
$$
where $\eta(x)$ denotes the unit normal at $x\in\Gamma$ pointing towards $\Omega_{n\ell}$ and $\partial_\eta u$ stands for derivativve in the direction of $\eta$ coming from the local part.
Since $\langle \partial E[u],\varphi\rangle=\partial_\varphi E(u)$, we arrive to a problem consisting of a local heat equation with a nonlocal source term in the \lq\lq local'' part of the  domain,
\begin{equation}  \label{eq:main.Cauchy.local}
\begin{cases}
\displaystyle u_t(x,t)=\Delta u +\alpha\int_{\Omega_{n\ell}} J(x-y)(u(y,t)-u(x,t))\,{\rm d}y,\qquad&x\in \Omega_\ell,\ t>0,
\\[6pt]
\displaystyle \partial_\eta u(x,t)=0,\qquad & x\in \Gamma, \ t>0,
\end{cases}
\end{equation}
together with a nonlocal heat equation in the \lq\lq nonlocal'' part of the domain, 
\begin{equation}  \label{eq:main.Cauchy.nonlocal}
\begin{aligned}
u_t(x,t)= \alpha\int_{\Omega_\ell} J(x-y)(u(y,t)-u(x,t))\, {\rm d}y + 2\alpha\int_{\Omega_{n\ell}} J(x-y) (u(y,t)-& u(x,t))\, {\rm d}y,
\\[6pt]
&
\quad x \in \Omega_{n\ell},\ t>0,
\end{aligned}
\end{equation}
plus an initial condition
$u(\cdot,0) = u_0 (\cdot)$ in $\mathbb{R}^N$.

From a probabilistic viewpoint (particle systems) in this model particles may jump (according with the probability density $J(x-y)$) when the initial point or the target point, $x$ or $y$,
belongs to the nonlocal region $\Omega_{n\ell}$, and also move according to Brownian motion (with a reflection at $\Gamma$) in the local region
$\Omega_\ell$. Notice that there is some interchange of mass between $\Omega_{n\ell}$ and $\Omega_\ell$ since particles may jump across $\Gamma$.

Notice that we do not impose any continuity to solutions of~\eqref{eq:main.Cauchy.local}--\eqref{eq:main.Cauchy.nonlocal} across the interface $\Gamma$ that separates the local and nonlocal domains.
In fact, solutions can be discontinuous across $\Gamma$ even if the initial condition is smooth.

 As precedents for our study we quote \cite{delia2,delia3,Du}. In \cite{delia2} local and nonlocal problems are coupled trough a prescribed region in which
both kind of equations overlap (the data from the nonlocal domain is used as a Dirichlet boundary condition for the local part and viceversa). This kind of coupling gives some continuity in the overlapping region but does not preserve the total mass.
In \cite{delia2} and \cite{Du} numerical schemes using local and nonlocal equations are developed and used in order to improve the computational
accuracy when approximating a purely nonlocal problem.

For this problem we have the following result:
\begin{Theorem}\label{thm:existence_uniqueness_modelo_1.intro}
	Given $u_0\in L^1(\mathbb{R}^N)$, there exists a unique $u\in C([0,\infty): L^1(\mathbb{R}^N))$ solving~\eqref{eq:main.Cauchy.local}--\eqref{eq:main.Cauchy.nonlocal} such that $u(\cdot,0)=u_0 (\cdot)$.
	The mass of the solution is conserved, $\int_{\mathbb{R}^N} u(\cdot,t)= \int_{\mathbb{R}^N} u_0 $. Moreover, a comparison
	principle holds: if $u_0 \geq v_0$ then the corresponding solutions verify $u \geq v$ in $\mathbb{R}^N \times \mathbb{R}_+$.
\end{Theorem}

\subsection{The Neumann problem}
Let us now present the version of this problem with boundary condition in dimension $N\geq 1$. Let us take a bounded smooth domain $\Omega\subset \mathbb{R}^N$ that is itself divided into two other subsets $\Omega_\ell$ and $\Omega_{n\ell}$ by a smooth hypersurface $\Gamma$. Again we can define an energy functional
\begin{equation}\label{def:funcional_neumann.intro}
E(u):=\int_{\Omega_\ell} \frac{|\nabla u|^2}{2} + \frac{\alpha}{2}\int_{\Omega_{n\ell}}\int_{\Omega} J(x-y)(u(y)-u(x))^2\,{\rm d}y{\rm d}x.
\end{equation}
Associated with this energy we obtain an evolution problem with a \lq\lq local'' part
\begin{equation}  
\label{eq:main.Neumann.local}
\begin{cases}
\displaystyle u_t(x,t)=\Delta u (x,t) + \alpha\int_{\Omega_{n\ell}} J(x-y)(u(y,t)-u(x,t))\,{\rm d}y,\ &x\in \Omega_\ell,\ t>0,
\\[6pt]
\displaystyle \partial_\eta u(x,t)=0,\qquad & x\in \partial \Omega_\ell, \ t>0,
\end{cases}
\end{equation}
and a \lq\lq nonlocal'' one, 
\begin{equation}  
\label{eq:main.Neumann.nonlocal}
\begin{aligned}
\displaystyle u_t(x,t)=\alpha\int_{\Omega_\ell} J(x-y)(u(y,t)-u(x,t))\, {\rm d}y
+ 2\alpha\int_{\Omega_{n\ell}}\!\! J(x-y)(u(y,t)-&u(x,t))\, {\rm d}y, \\[6pt]
&
\quad x \in \Omega_{n\ell},\ t>0,
\end{aligned}
\end{equation}
plus an initial condition $u(\cdot,0) = u_0(\cdot)$ in $\Omega$.

Notice that in this model there are no individuals that may jump into $\Omega$ coming from the outside side $\mathbb{R}^N \setminus \Omega$ nor individuals that jump
from $\Omega$ into the exterior side $\mathbb{R}^N \setminus \Omega$. It is in this sense that we call~\eqref{eq:main.Neumann.local}--\eqref{eq:main.Neumann.nonlocal} a Neumann
type problem.

For this problem we also have existence and uniqueness  of  solutions and a comparison principle. Moreover, as it is expected for Neumann
boundary conditions, we also have conservation of mass.

\begin{Theorem}\label{thm:existence_uniqueness_neumann.intro}
	Given $u_0\in L^1(\Omega)$, there exists a unique $u\in C([0,\infty): L^1(\Omega))$ solving~\eqref{eq:main.Neumann.local}--\eqref{eq:main.Neumann.nonlocal} such that $u(\cdot,0)=u_0 (\cdot)$. This solution conserves mass.
Moreover, a comparison principle holds.
\end{Theorem}

\subsection{The Dirichlet problem}
As for the Dirichlet case let us take a bounded smooth domain $\Omega\subset \mathbb{R}^N$ that is itself divided into two subsets $\Omega_\ell$ and $\Omega_{n\ell}$ by a smooth hypersurface $\Gamma$.
The Dirichlet version of the functional reads as
\begin{equation}\label{def:funcional_dirichlet.intro}
E(u):=\int_{\Omega_\ell} \frac{|\nabla u|^2}{2} + \frac{\alpha}{2}\int_{\Omega_{n\ell}}\int_{\mathbb{R}^N} J(x-y)(u(y)-u(x))^2\, {\rm d}y{\rm d}x,
\end{equation}
extending $u$ by zero outside $\Omega$ (and hence also on $\partial \Omega$). Notice that in the nonlocal part we have integrated $y$ in the whole $\mathbb{R}^N$.
The associated evolution problem again has a local part,
\begin{equation}  \label{eq:main.Dirichlet.local}
\begin{cases}
\displaystyle u_t(x,t)=\Delta u (x,t) + \alpha\int_{\mathbb{R}^N\setminus \Omega_\ell} J(x-y)(u(y,t)-u(x,t))\,{\rm d} y, &x\in \Omega_\ell,\, t>0,
\\[6pt]
\displaystyle \partial_\eta u(x,t)=0,\qquad & x\in \Gamma, \ t>0,
\end{cases}
\end{equation}
plus a nonlocal one
\begin{equation}  \label{eq:main.Dirichlet.nonlocal}
\begin{aligned}
\displaystyle u_t(x,t)=\alpha\int_{\mathbb{R}^N\setminus \Omega_{n\ell}} J(x-y)(u(y,t)-u(x,t))\, {\rm d}y
+ 2\alpha\int_{\Omega_{n\ell}} J(x-y) & (u(y,t)-u(x,t))\, {\rm d}y,
\\[6pt]
& \quad x \in \Omega_{n\ell},\ t>0,
\end{aligned}
\end{equation}
plus the Dirichlet \lq\lq boundary'' condition
\begin{equation}  \label{eq:main.Dirichlet.boundary}
\displaystyle u(x,t)= 0,\quad  x\in \mathbb{R}^N\setminus \Omega,\  t>0,
\end{equation}
and the initial condition $u(\cdot,0)=u_0 (\cdot)$ in $\Omega$. 

In this model we have that individuals may jump outside $\Omega$ but they instantaneously die there since we have that the density
$u$ vanishes identically in $(\mathbb{R}^N\setminus \Omega)\times \mathbb{R}_+$.

For this problem we also have existence and uniqueness of solutions and a comparison principle, but, of course, there is no
conservation of mass.
\begin{Theorem}\label{thm:existence_uniqueness_Dirichlet.intro}
Given $u_0\in L^1(\Omega)$, there exists a unique $u\in C([0,\infty): L^1(\Omega))$ solving~\eqref{eq:main.Dirichlet.local}--\eqref{eq:main.Dirichlet.boundary} such that $u(\cdot,0)=u_0(\cdot)$. Moreover, a comparison principle holds.
\end{Theorem}

\subsection{Asymptotic behavior.} It is well known that solutions of the local heat equation~\eqref{heat.equation} have a polynomial time decay or the Cauchy problem  and an exponential decay
(to the mean value of the initial condition or to zero) for the Neumann and the Dirichlet problems. The same is true for solutions of the nonlocal heat equation~\eqref{eq:nonlocal.equation}, though in the case of the Cauchy problem we have to ask the second moment of the kernel
\begin{equation*}
\label{eq:def.second.moment}
M_2(J):=\int_{\mathbb{R}^N}J(z)|z|^2\,{\rm d}z
\end{equation*}
to be finite as in~\cite{Chasseigne-Chaves-Rossi-2006}. Our local/nonlocal model reproduces these behaviours.
\begin{Theorem}\label{thm:asymp.intro}
{\rm (a)} Let $u$ be a solution of the Cauchy problem~\eqref{eq:main.Cauchy.local}--\eqref{eq:main.Cauchy.nonlocal} in $\mathbb{R}^N$ with integrable initial data. If $M_2(J)<\infty$,  for any $p\in(1,\infty)$ there is a constant $C$ such that 
\begin{equation} \label{decay.cauchy}
\|u(\cdot, t ) \|_{L^{p} (\mathbb{R}^N)} \leq C t^{-\frac{N}{2} (1 - \frac{1}{p})}.
\end{equation}

\noindent{\rm (b)} Let $u$ be a solution of the Dirichlet problem problem~\eqref{eq:main.Dirichlet.local}--\eqref{eq:main.Dirichlet.boundary} in an $N$-dimensional domain with integrable initial data. For any $p\in[1,\infty)$ there are positive constants $C$ and $\lambda$ such that 
\begin{equation} \label{decay.Dirichlet}
\|u(\cdot, t ) \|_{L^{p} (\Omega)} \leq C e^{-\lambda t}.
\end{equation}

\noindent{\rm (c)} Let $u$ be a solution of the Neumann problem~\eqref{eq:main.Neumann.local}--\eqref{eq:main.Neumann.nonlocal} in an $N$-dimensional domain. For any $p\in[1,\infty)$ there are positive constants $C$ and $\beta>0$ such that 
\begin{equation} \label{decay}
\left\|u(\cdot, t ) -  |\Omega|^{-1}\int_{\Omega} u_0 \right\|_{L^{p} (\Omega)} \leq C e^{-\beta t}.
\end{equation}
\end{Theorem}

\subsection{Rescaling the kernel.} Our aim is to recover the usual local problems from our nonlocal ones when
we rescale the kernel according to
\begin{equation}
\label{eq:definition.J.epsilon}
J^\varepsilon(z)=\varepsilon^{-(N+2)}J(z/\varepsilon).
\end{equation}
It is at this point where we choose the constant $\alpha$ that appears in front of our nonlocal terms as
$$
\alpha=1/M_2(J).
$$
Now, for a fixed initial condition $u_0$ and for each $\varepsilon>0$ our evolution problems (Cauchy, Neumann or Dirichlet) have a solution.
Our goal is to look for the limit as $\varepsilon \to 0$ of these solutions to recover in this limit procedure the local heat equation.
Notice that as $\varepsilon$ becomes small the support of the kernel $J^\varepsilon$ shrinks, hence the non locality of the operator
becomes weaker as $\varepsilon$ becomes smaller.
As precedents where this kind of limit procedure is performed we quote
\cite{BLGneu,BLGdir,ElLibro,BChQ,Canizo-Molino-2018,Chass,CER,Cortazar-Elgueta-Rossi-Wolanski,Andres}.
One of the main difficulties here is that we do not have any continuity of the solutions to our nonlocal
equations across the interface that separates the local and nonlocal domains, while the expected limit is smooth
across the interface (being a solution to the heat equation in the whole domain).

\begin{Theorem}\label{thm:rescales.intro} Let $u_0\in L^2(\Omega)$ (with $\Omega=\mathbb{R}^N$ in the case of the Cauchy problem). For each $\varepsilon>0$, let $u^\varepsilon$ be the solution of  any of the three previously mentioned problems, Cauchy, Neumann, or Dirichlet with initial data $u_0$. Then, as $\varepsilon \to 0$,
\begin{equation} \label{conver.intro}
u^\varepsilon \to u \quad \mbox{strongly in } L^2,
\end{equation}
where $u$ is the solution to the corresponding problem (Cauchy, Neumann, or Dirichlet, in the two last cases with zero boundary condition) for the local heat equation~\eqref{heat.equation} in $\Omega\times\mathbb{R}_+$ with the same initial condition. 
\end{Theorem}

\medskip

The rest of the paper is organized as follows: in Section \ref{sect-prelim} we collect some preliminary results
and prove an inequality that will be the key in our arguments;
in the following sections we prove our main results concerning existence, uniqueness and properties of the model in its three versions
(Cauchy, Neumann, Dirichlet).
We gather the results according to the problem we deal with and hence
in Section \ref{sect-Cauchy} we study the Cauchy problem (including its asymptotic behaviour and the limit when we rescale the kernel); in Section \ref{sect-Neumann} the Neumann problem and
finally in Section \ref{sect-Dirichlet} the Dirichlet problem.

\section{Preliminaries} \label{sect-prelim}
\setcounter{equation}{0}
First, we present a very useful result that we state in its more general form. This result says that we can control
the purely nonlocal energy by our local/nonlocal one.

\begin{Lemma}\label{lemma:preliminar}
Let $\Omega\subset\mathbb{R}^N$ be a smooth domain, $\Omega_\ell\subset \Omega$ a smooth convex subdomain and $\Omega_{n\ell}=\Omega\setminus \overline{\Omega_\ell}$. Let $u\in\mathcal{H}:=\{u\in L^2(\Omega): u|_{\Omega_{\ell}}\in H^1(\Omega_{\ell}) \}$. Then, for any $k\in (0,1/M_2(J))$, 
\begin{equation}
\label{eq:fundamental.inequality}
\begin{array}{l}
\displaystyle
\int_{\Omega_\ell} |\nabla u|^2 + \frac{1}{2M_2(J)}\int_{\Omega_{n\ell}}\int_{\Omega} J(x-y)(u(y)-u(x))^2\, {\rm d}y{\rm d}x
\\ [4mm]
\qquad \geq \displaystyle k\int_{\Omega}\int_{\Omega} J(x-y)(u(y)-u(x))^2\, {\rm d}y{\rm d}x.
\end{array}
\end{equation}
\end{Lemma}

\begin{proof}
Thanks to the symmetry of the kernel $J$, our inequality is equivalent to showing that
\begin{equation}
\int_{\Omega_\ell} |\nabla u|^2  -k\int_{\Omega_\ell}\int_{\Omega_\ell} J(x-y)(u(y)-u(x))^2\, {\rm d}y {\rm d}x\geq 0,
\end{equation}
if $k$ is small enough.
This is important because we stick to the domain $\Omega_\ell$, where  $u$ belongs to $H^1$, so we can express it as the integral of a derivative and make computations with it.

After a change of variables, using Jensen's inequality we get
\begin{equation}
\begin{aligned}
\displaystyle\int_{\Omega_\ell}\int_{\Omega_\ell} J(x-y)&(u(y)-u(x))^2\,{\rm d}y{\rm d}x = \int_{\Omega_\ell}\int_{\Omega_\ell-x} J(z)(u(x+z)-u(x))^2\,{\rm d}z{\rm d}x
\\[6pt]
\displaystyle &=\int_{\Omega_\ell}\int_{\Omega_\ell-x} J(z)\left(\int_0^{1} \nabla u(x+sz)\cdot z\, {\rm d}s\right)^2\, {\rm d}z{\rm d}x\\[6pt]
\displaystyle&\leq \int_{\Omega_\ell}\int_{\Omega_\ell-x}\int_0^{1} J(z)|z|^2 |\nabla u(x+sz)|^2\, {\rm d}s{\rm d}z{\rm d}x=:\mathcal{A}.
\end{aligned}
\end{equation}
If we define the sets $Z:=\{z\in\mathbb{R}^N:\text{there exists }x\in \Omega_\ell \text{ such that } x+z\in\Omega_\ell\}$ and $\Omega_{z,s}:=\{x\in\Omega_\ell:x+sz \in\Omega_\ell\}$ then we can apply Fubini to see that 
$$
\mathcal{A}=\int_{Z} \int_{0}^1\int_{\Omega_{z,s}} J(z)|z|^2 |\nabla u(x+sz)|^2\, {\rm d}x{\rm d}s{\rm d}z.
$$
Now we define, for fixed $z,s$ the variable $w:=x+sz$, which means that the set $\Omega_{z,s}$ can be described as $W_{z,s}:=\{w\in\Omega_\ell : \text{ there exists } x\in\Omega_\ell \text{ such that } w=x+sz \}\subset \Omega_\ell$. Hence,
$$
\begin{aligned}
\mathcal{A}=\int_{Z} \int_{0}^1\int_{W_{z,s}} &J(z)|z|^2 |\nabla u(w)|^2\,{\rm d}w{\rm d}s{\rm d}z \leq\int_{Z} \int_{0}^1\int_{\Omega_\ell} J(z)|z|^2 |\nabla u(w)|^2\,{\rm d}w{\rm d}s{\rm d}z
\\[6pt]
\displaystyle &\leq\int_{\mathbb{R}^N} \int_{0}^1\int_{\Omega_\ell} J(z)|z|^2 |\nabla u(w)|^2\,{\rm d}w{\rm d}s{\rm d}z = M_2(J)\int_{\Omega_\ell} |\nabla u|^2.
\end{aligned}
$$
The result follows taking $k$ small enough.
\end{proof}

It is worth noting that estimate~\eqref{eq:fundamental.inequality} scales well with the rescaled version of the kernel $J^\varepsilon$ given by~\eqref{eq:definition.J.epsilon}, 
since $M_2(J^\varepsilon)=M_2(J)$. Hence we can take the same constant $k$ for all $\varepsilon$.  This will be helpful when studying the asymptotic behaviour of the solutions of these problems and also for the convergence of these problems to the corresponding local ones.

The following lemma will be needed later on to study the limit behaviour under rescales of the kernel for the three different problems and is an adaptation of the results that can be found in~\cite{BLGneu,BLGdir, ElLibro}. In the following we will denote by $\bar{f}$ the extension by zero of a function  $f$ outside  our domain $\Omega$.

\begin{Lemma}\label{lemma:basic_convergence_lemma_for_rescaled}
Let $\Omega$ be a bounded domain and $\{f^\varepsilon\}$ a sequence of functions in $L^2(\Omega)$ such that
$$
\displaystyle \int_{\Omega}\int_{\Omega} J^\varepsilon(x-y)(f^\varepsilon(y)-f^\varepsilon(x))^2\, {\rm d}y{\rm d}x \leq C
$$
for a positive constant $C$ and $\{f^\varepsilon\}$ is weakly convergent in $L^2(\Omega)$ to $f$ as $\varepsilon$ goes to 0.
Then $$|\nabla f|\in L^2(\Omega)$$ and moreover
$$
\lim\limits_{\varepsilon\to 0} \left( J^{1/2}(z) \chi_{\Omega}(x+\varepsilon z) \frac{\bar{f}^\varepsilon(x+\varepsilon z) - f^\varepsilon(x)}{\varepsilon} \right) = J^{1/2}(z)z\nabla f(x)
$$
weakly in $L^2_x(\Omega)\times L^2_z(\mathbb{R}^N)$.
\end{Lemma}
\begin{proof}
Changing variables $y=x+\varepsilon z$ we obtain
$$
\displaystyle \int_{\Omega}\int_{\mathbb{R}^N} J(z) \chi_{\Omega}(x+\varepsilon z) \frac{(\bar{f}^\varepsilon(x+\varepsilon z,t)-f^\varepsilon(x,t))^2}{\varepsilon^2}\,{\rm d}z{\rm d}x \leq C,
$$
which already provides, for a certain function $h=h(x,z)$, the stated weak convergence to a $J^{1/2}(z)h(x,z)$. Having this weak convergence we can multiply the quantity
$$
J^{1/2}(z) \chi_{\Omega}(x+\varepsilon z) \frac{\bar{f}^\varepsilon(x+\varepsilon z) - f^\varepsilon(x)}{\varepsilon}
$$
by two test functions $\varphi(x)\in C_c^\infty(\Omega)\subset L^2(\mathbb{R}^N)$ and $\psi(z)\in C_c^\infty(\mathbb{R}^N)\subset L^2(\mathbb{R}^N)$, integrate and pass to the limit to obtain that
\begin{equation}\label{eq:lema_preliminar_1}
\begin{aligned}
\displaystyle \int_{\mathbb{R}^N} &J^{1/2}(z)\psi(z) \int_{\Omega}  \chi_{\Omega}(x+\varepsilon z) \frac{\bar{f}^\varepsilon(x+\varepsilon z) - f^\varepsilon(x)}{\varepsilon} \varphi(x) \,{\rm d}x{\rm d}z
\\[6pt]
\displaystyle &\to \int_{\mathbb{R}^N} J^{1/2}(z)\psi(z) \int_{\Omega} h(x,z)\varphi(x)\,{\rm d}x{\rm d}z,
\end{aligned}
\end{equation}
and it is now when we note that since $\psi(z)$ has compact support the integral over $\mathbb{R}^N$ in $z$ is really an integral over a compact set, so there exists $\varepsilon$ small enough such that
$$
J^{1/2}(z) \chi_{\Omega}(x+\varepsilon z) = J^{1/2}(z)
$$
for all $z\in \mathbb{R}^N$ and all $x\in \mathop{\rm supp}(\varphi)$. Then
$$
\begin{aligned}
\displaystyle \int_{\mathbb{R}^N} &J^{1/2}(z)\psi(z) \int_{\Omega}  \chi_{S}(x+\varepsilon z) \frac{\bar{f}^\varepsilon(x+\varepsilon z) - f^\varepsilon(x)}{\varepsilon} \varphi(x) \,{\rm d}x{\rm d}z
\\[6pt]
&= \displaystyle \int_{\mathbb{R}^N} J^{1/2}(z)\psi(z) \int_{\mathop{\rm supp}(\varphi)}\frac{\bar{f}^\varepsilon(x+\varepsilon z) - f^\varepsilon(x)}{\varepsilon} \varphi(x) \,{\rm d}x{\rm d}z
\\[6pt]
&= \displaystyle - \int_{\mathbb{R}^N} J^{1/2}(z)\psi(z) \int_{\Omega}f^\varepsilon(x)\frac{\varphi(x) - \bar{\varphi}(x-\varepsilon z)}{\varepsilon} \,{\rm d}x{\rm d}z.
\end{aligned}
$$
Using~\eqref{eq:lema_preliminar_1} and the fact that $\{f^\varepsilon\}$ converges weakly to $f$ in $L^2(S)$ we have that
$$
\displaystyle - \int_{\mathbb{R}^N} J^{1/2}(z)\psi(z) \int_{\Omega}f\nabla\varphi(x) \,{\rm d}x{\rm d}z = \displaystyle \int_{\mathbb{R}^N} J^{1/2}(z)\psi(z) \int_{\Omega} h(x,z)\varphi(x)\,{\rm d}x{\rm d}z
$$
and from this point it is easy to conclude what remains of the lemma.
\end{proof}

\section{The Cauchy problem} \label{sect-Cauchy}
\setcounter{equation}{0}

In this section, we will prove existence and uniqueness of solutions to the Cauchy problem~\eqref{eq:main.Cauchy.local}--\eqref{eq:main.Cauchy.nonlocal} with initial data $u_0\in L^1(\mathbb{R}^N)$, conservation of mass, the asymptotic polynomial decay of the $L^p$ norms in time, and the convergence of this problem to the local one when we rescale the kernel.

%

\subsection{Existence and uniqueness}
To prove existence and uniqueness of a solution the idea is to use a fixed point argument as follows: Given a function $v$ defined for $\Omega_{\ell}$ we solve for $\Omega_{n\ell}$, in a function we will call $z$. With the obtained function we solve back for $\Omega_{\ell}$ in a function $w$. This can be regarded as an operator $T$ that satisfies $T(v)=w$. This is the operator for which we will look for a fixed point via contraction in adequate norms (we will use this technique several times), meaning that there must exist a $v=T(v)$ solving the equation for $\Omega_{\ell}$ with its corresponding $z$ solving the equation for $\Omega_{n\ell}$.

We will write this argument for an initial condition $u_0\in L^1(\mathbb{R}^N)$. Let us define, for a fixed finite $t_0>0$ the norms
$$
\|v\|_\ell= \sup\limits_{t\in [0,t_0]} \|v(\cdot, t)\|_{L^1(\Omega_\ell)} \qquad  \text{and}
\qquad \|z\|_{n\ell} = \sup\limits_{t\in [0,t_0]} \|z(\cdot, t)\|_{L^1(\Omega_{n\ell})}.
$$

Given $t_0>0$ to be chosen later, we define an operator $T_1:L^1(\Omega_\ell\times (0,t_0)))\to L^1(\Omega_{n\ell}\times(0,t_0))$ as $T_1(v)=z$, where $z$ is the unique solution to
\begin{equation}\label{eq:solucion_z}
\begin{cases}
\displaystyle z_t(x,t)=\alpha\int_{\Omega_\ell} J(x-y)(v(y,t)-z(x,t))\,{\rm d}y
+ 2\alpha \int_{\Omega_{n \ell}} J(x-y)(z(y,t)-z(x,t))\,{\rm d}y,
\\[6pt]
\displaystyle \hskip10cm\qquad x\in\Omega_{n \ell},\ t\in(0,t_0),
\\[6pt]
\displaystyle z(x,0) = u_0(x),\qquad x\in \Omega_{n \ell} .
\end{cases}
\end{equation}
Let us check that this problem has indeed a unique solution. In addition, we will study its dependence on the data $v$.
\begin{Lemma} Let $t_0\in (0,1/(5\alpha))$.  
Given $v\in L^1(\Omega_\ell\times (0,t_0)))$ there exists a unique $z\in L^1(\Omega_{n\ell}\times(0,t_0))$ that solves~\eqref{eq:solucion_z} for some $t_0$ small enough. Moreover, if $z_1$ and $z_2$ are the solutions corresponding respectively to $v_1$ and $z_2$,  then
\begin{equation}
\label{eq:contraction.z.v}
\|z_1-z_2\|_{n\ell} \leq \frac{\alpha t_0}{1-5\alpha t_0} \|v_1-v_2\|_\ell.
\end{equation}
\end{Lemma}
\begin{proof}
To show existence and uniqueness we use a fixed point argument. We define an operator $S_v:L^1(\Omega_{n\ell}\times (0,t_0)))\to L^1(\Omega_{n\ell}\times(0,t_0))$ through
$$
\begin{array}{ll} \displaystyle
S_v(z)(x,t):=u_0(x) + \alpha\int_0^t \int_{\Omega_\ell} J(x-y)(v(y,s)-z(x,s))\,{\rm d}y{\rm d}s
\\ [4mm]
\displaystyle\qquad \qquad \qquad + 2\alpha\int_0^t \int_{\Omega_{n\ell}} J(x-y)(z(y,s)-z(x,s))\,{\rm d}y{\rm d}s \quad\text{for }x\in \Omega_{n\ell},\, t\in(0,t_0).
\end{array}
$$
An easy computation shows that
$$
\begin{array}{ll}
\displaystyle
\|S_v(z_1)-S_v(z_2)\|_{n\ell} \leq \alpha \sup\limits_{t\in[0,t_0]} \int_{\Omega_{n\ell}}
\int_0^t \int_{\Omega_\ell} J(x-y)|z_2(x,s)-z_1(x,s)|\,{\rm d}y{\rm d}s{\rm d}x
\\ [4mm]
\displaystyle\qquad\qquad\qquad\qquad\qquad  \qquad+ 2\alpha \sup\limits_{t\in[0,t_0]} \int_{\Omega_{n\ell}} \int_0^t \int_{\Omega_{n\ell}}
J(x-y)|z_1(y,s)-z_2(y,s)|\,{\rm d}y{\rm d}s{\rm d}x
\\ [4mm]
\displaystyle\qquad \qquad\qquad\qquad\qquad\qquad + 2\alpha\sup\limits_{t\in[0,t_0]} \int_{\Omega_{n\ell}}\int_0^t
\int_{\Omega_{n\ell}} J(x-y)|z_2(x,s)-z_1(x,s)|\,{\rm d}y{\rm d}s{\rm d}x,
\end{array}
$$
but here we recall that the integral of the kernel $J$ is always lesser or equal than 1, and apply Fubini's theorem to obtain
$$
\|S_v(z_1)-S_v(z_2)\|_{n\ell} \leq 5\alpha t_0\|z_1-z_2\|_{n\ell}.
$$
Choosing $t_0\leq 1/(5\alpha)$, $S_v$ is a strict contraction, and hence has a unique fixed point. 

As for the dependence on the data, since $z_1=S_{v_1}(z_1)$ and $z_2=S_{v_2}(z_2)$, a computation similar to the one we have just performed gives
$$
\|z_1-z_2\|_{n\ell} \leq 5\alpha t_0\|z_1-z_2\|_{n\ell} + \alpha t_0\|(v_1-v_2)(x,t)\mathbb{R}
$$
which yields~\eqref{eq:contraction.z.v}.
\end{proof}

Now it is time to go back to $\Omega_\ell$.  We define 
$T_2:L^1(\Omega_{n\ell}\times (0,t_0)))\to L^1(\Omega_{\ell}\times(0,t_0))$ as $T_2(z)=w$, where $w$ is the unique solution to
\begin{equation}\label{eq:solucion_w}
\begin{cases}
\displaystyle w_t(x,t)= \Delta w (x,t) -C_1w(x,t)A(x)+C_2 B_z(x,t),\quad&(x,t)\in \Omega_{\ell}\times(0,t_0),
\\[6pt]
\displaystyle -\partial_\eta w(x,t)=0,\quad & x\in \Gamma, \, t\in(0,t_0),
\\[6pt]
\displaystyle w(x,0) = u_0(x),\qquad &x\in\Omega_\ell,
\end{cases}
\end{equation}
with
$$
A(x):=\int_{\Omega_{n\ell}} J(x-y)\, {\rm d}y\quad\text{and }B_z(x,t) := \int_{\Omega_{n\ell}} J(x-y)z(y,t)\,{\rm d}y.
$$

\begin{Lemma} Let $t_0\in (0,1/(2\alpha))$.  
Given $z\in L^1(\Omega_{n\ell}\times (0,t_0)))$ there exists a unique $w\in L^1(\Omega_{\ell}\times(0,t_0))$ that solves~\eqref{eq:solucion_z} for some $t_0$ small enough. Moreover, if $w_1$ and $w_2$ are the solutions corresponding respectively to $z_1$ and $z_2$,  then
\begin{equation}
\label{eq:contraction.z.v}
\|w_1-w_2\|_{n\ell} \leq \frac{\alpha t_0}{1-\alpha t_0} \|z_1-z_2\|_\ell.
\end{equation}
\end{Lemma}

\begin{proof}
Existence and uniqueness of solutions is well known, see \cite{Evans}. The contraction property follows in a similar way as before, taking into account the condition at the boundary and the source term. This time we obtain the estimate
$$
\|w_1-w_2\|_\ell\leq \frac{\alpha t_0}{1-\alpha t_0}\|z_1-z_2\|_{n\ell},
$$
which is a contraction given $t_0\leq 1/(2\alpha)$.
\end{proof}

Thus, combining the two previous lemmas, we have obtained the following theorem.
\begin{Theorem}\label{thm:existence_uniqueness_modelo_1}
Given $u_0\in L^1(\mathbb{R}^N)$,  there exists a unique solution to problem~\eqref{eq:main.Cauchy.local}--\eqref{eq:main.Cauchy.nonlocal} which has $u_0$ as initial datum.
\end{Theorem}

\begin{proof}
First, we keep $t_0\leq 1/(6\alpha)$). 
If we compose now the two operators
$$w=T(v):=T_2(T_1(v))$$
is easy to obtain
$$
\|w_1-w_2\|_\ell \leq \frac{\alpha^2t_0^2}{(1-\alpha t_0)(1-5\alpha t_0)}\|v_1-v_2\|_\ell,
$$
which again is contraction given $t_0\leq 1/(6\alpha)$. Therefore, there is a fixed point
that gives us a unique solution in $(0,t_0)$. Now, using that the fixed point argument can be
iterated we obtain a global solution for our problem.
\end{proof}

There is also an alternative approach to prove existence of solutions for this problem. Applying the linear semi-group theory, see~\cite{ElLibro}, we can define the operator
$$
B_J(u)= \begin{cases}
\displaystyle
             -\Delta u (x) - \alpha\int_{\Omega{n\ell}} J(x-y)(u(y)-u(x))\,{\rm d}y  & \text{if }x\in\Omega_\ell, \\[8pt]
             \displaystyle  -2\alpha\int_{\Omega_{n\ell}} J(x-y)(u(y)-u(x))\,{\rm d}y - \alpha\int_{\Omega_\ell} J(x-y)(u(y)-u(x))\,{\rm d}y  & \text{if }x\in\Omega_{n\ell}
             \end{cases}
$$
with the constraint that $\partial_\eta u=0$ on $\Gamma$ and let $D(B_J)$ denote its domain. Following~\cite{ElLibro} we will see that this operator is completely accretive and satisfies the range condition $L^2(\mathbb{R}^N)\in R(I+B_J)$. This will imply that for any $\phi\in L^2(\mathbb{R}^N)$ there exists a $u\in D(B_J)$ such that $u+B_J(u)=\phi$ and the resolvent $(I+B_J)^{-1}$ is a contraction in $L^p(\mathbb{R}^N)$ for every $1\leq p\leq \infty$. After that, Crandall-Ligget's Theorem and the linear semi-group theory will give existence and uniqueness of a \textit{mild} solution of our evolution problem. In what follows we will use notations form semi-group theory, see~\cite{ElLibro}.

\begin{Theorem}
The operator $B_J(u)$ is completely accretive and satisfies the range condition $$L^2(\mathbb{R}^N)\in R(I+B_J).$$
\end{Theorem}

\begin{proof}
To show that the operator is completely accretive it is enough to see that for every given $u_i\in D(B_J)$, $i=1,2$, and $q\in C^\infty(\mathbb{R})$ such that $0\leq q'\leq 1$, $\mathop{\rm supp}(q)$ is compact and $0\not\in \mathop{\rm supp}(q)$ we have that
$$
\int_{\mathbb{R}^N}\left(B_J (u_1(x)) - B_J (u_2(x))\right)\cdot q(u_1(x)-u_2(x))\,{\rm d}x\geq 0.
$$
To see this is not difficult through a change of variables $x \leftrightarrow y$, Fubini and the Mean Value Theorem that give us
$$
q(u_1(x)-u_2(x)) =  q' (\xi)\cdot (u_1(x)-u_2(x))
$$
for some real intermediate real number $\xi$.

To show the range condition let us take first $\phi \in L^\infty(\mathbb{R}^N)$ and define the auxiliary operator
$$
A_{n,m}(u):= T_c(u) + B_J(u)  + \frac{1}{n}u^+ - \frac{1}{m}u^-
$$
where $T_c(u):=\min(c, \max(u,-c))$ is the function $u$ truncated between $-c$ and $c$. This operator is continuous monotone, and more importantly it is easy to check that it is coercive in $L^2(\mathbb{R}^N)$. Then, by~\cite{Brezis-1968}, there exists a $u_{n,m}\in L^2(\mathbb{R}^N)$ such that
$$
T_c(u_{n,m}) + B_J(u_{n,m})  + \frac{1}{n}u_{n,m}^+ - \frac{1}{m}u_{n,m}^-=\phi.
$$
Let us also define the following relation. We will write $f\ll g$ if and only if
$$
\int_{\mathbb{R}^N} j(f) \leq \int_{\mathbb{R}^N} j(g)
$$
for every $j:\mathbb{R}\to [0,\infty]$ convex, lower semi-continuous and with $j(0)=0$.

Using the monotonicity of $$ B_J(u_{n,m})  + \frac{1}{n}u_{n,m}^+ - \frac{1}{m}u_{n,m}^-$$ we have that $T_c(u_{n,m})\ll \phi$, so taking $c>\|\phi\|_{L^\infty(\mathbb{R}^N)}$ we see that $u_{n,m}\ll \phi$ and
$$
u_{n,m} + B_J(u_{n,m})  + \frac{1}{n}u_{n,m}^+ - \frac{1}{m}u_{n,m}^-=\phi.
$$

Now we will see that $u_{n,m}$ is non-decreasing in $n$ and non-increasing in  $m$ in order to pass to the limits. We will show the ideas for the monotonicity in $n$, since is similar for $m$. We define $w:=u_{n,m} - u_{n+1,m}$ and this $w$ satisfies
$$
w+B_J(w)+ \frac{1}{n}u_{n,m}^+ - \frac{1}{n+1}u_{n+1,m}^+ + \frac{1}{m}u_{n+1,m}^- - \frac{1}{m}u_{n,m}^- =0.
$$
We can now multiply by $w^+$ and integrate to obtain
$$
\int_{\mathbb{R}^N}(w^+)^2+\int_{\mathbb{R}^N}B_J(w)w^+\leq 0.
$$
Through already mentioned techniques is easy to check that the second integral is positive, meaning that necessarily $w^+=0$, meaning that $u_{n,m} \leq u_{n+1,m}$. Since for the parameter $m$ is similar we have the mentioned monotonicity. Thus, using that $u_{n,m}\ll \phi$, this monotonicity and monotone convergence for the term $B_J(u_{n,m})$, we pass to the limit $n\to \infty$ to obtain
$$
u_{m} + B_J(u_{m}) - \frac{1}{m}u_{m}^-=\phi
$$
and $u_m\ll\phi$. Passing again to the limit in $m$ we obtain
$$
u + B_J(u)=\phi.
$$

Now let $\phi\in L^2(\mathbb{R}^N)$ and $\phi_n$ a sequence in $L^\infty(\mathbb{R}^N)$ such that $\phi_n\to\phi$ in $L^2(\mathbb{R}^N)$. Then we have existence for $u_n=(I+B_J)^{-1}\phi_n$ by the previous steps and due to the complete accretiveness of the operator $u_n\to u$ in $L^2(\mathbb{R}^N)$ and $B_J(u_n)\to B_J(u)$ in $L^2(\mathbb{R}^N)$ (since the dual of $L^2$ is itself). We conclude then that $u+B_J(u)=\phi_n$.
\end{proof}

\subsection{Conservation of mass}
As expected, this model preserves the total mass of the solution. Formally, we have
$$
\begin{aligned}
\partial_t \int_{\mathbb{R}^N}u(x,t)\,{\rm d}x =&  \int_{\Omega_\ell} \Delta u (x,t)\,{\rm d}x + \alpha\int_{\Omega_\ell} \int_{\Omega_{n\ell}} J(x-y)(u(y,t)-u(x,t))\,{\rm d}y{\rm d}x
\\[6pt]
&+ \alpha\int_{\Omega_{n\ell}} \int_{\Omega_\ell} J(x-y)(u(y,t)-u(x,t))\,{\rm d}y{\rm d}x
\\[6pt]
 &+
  2  \alpha\int_{\Omega_{n\ell}} \int_{\Omega_{n\ell}} J(x-y)(u(y,t)-u(x,t))\, {\rm d}y{\rm d}x=0.
\end{aligned}
$$
The first integral is 0 thanks to the boundary condition on $u_x(0,t)$, the second and third ones add up to 0, changing variables $x$ and $y$ and using Fubini, and the last integral is equal to 0 due to the domain of integration and the symmetry of the kernel $J$. With all this and multiplying by a suitable test function our solution it is easy to prove the following theorem.

\begin{Theorem} \label{teo.conserva.masa.Cauchy}
The solution $u$ of problem~\eqref{eq:main.Cauchy.local}--\eqref{eq:main.Cauchy.nonlocal} with initial data $u_0\in L^1(\mathbb{R}^N)$ satisfies
$$
\int_{\mathbb{R}^N} u(x,t) \, {\rm d}x =\int_{\mathbb{R}^N} u_0(x) \, {\rm d}x\qquad \text{for all } t\geq 0.
$$
\end{Theorem}

\subsection{Comparison principle}

If we have two different solutions of the Cauchy problem problem~\eqref{eq:main.Cauchy.local}--\eqref{eq:main.Cauchy.nonlocal} then thanks to the linearity of the operator the difference between them is also a solution. It is also easy to see that given a non-negative initial data $u_0$ the solution keeps this non-negativity
(this follows from the fixed point construction of the solution or from the accretivity of the associated operator).
With this in mind, we state the following result.

\begin{Theorem}
If $u_0\geq 0$ then $u\geq 0$ in $\mathbb{R}^N\times\mathbb{R}_+$. Moreover, given two initial data $u_0$ and $v_0$, if $u_0\geq v_0$ then
$u\geq v$ in $\mathbb{R}^N\times\mathbb{R}_+$.
\end{Theorem}

\subsection{Asymptotic decay}

To study the decay of this problem we need a result that can be found in~\cite{Canizo-Molino-2018}.

\begin{Proposition}[\cite{Canizo-Molino-2018}]\label{proposition:molino}
Take the energy functional
$$
D_p^J(u)=\frac{p}{2}\int_{\mathbb{R}^N}\int_{\mathbb{R}^N}J(x-y)(u(y,t)-u(x,t))(|u|^{p-2}u(y,t)-|u|^{p-2}u(x,t))\,{\rm d}y{\rm d}x.
$$
Then, for every $u\in L^1\cap L^{\infty}(\mathbb{R}^N)$ and $p\geq 2$ there exists a positive constant $C$ such that
$$
D_p^J(u)\geq C\min\{\|u\|_{L^1(\mathbb{R}^N)}^{-p\gamma} \|u\|_{L^p(\mathbb{R}^N)}^{p(1+\gamma)}, \|u\|_{L^p(\mathbb{R}^N)}^p\}
$$
where $\gamma=2/N(p-1)$, and this bound provides a decay of the solutions of the evolution problem $$u_t=-(D_p^J)'(u)$$ of the form
$$
\|u(\cdot,t)\|_{L^p(\mathbb{R}^N)}^p \leq Ct^{-\frac{N(p-1)}{2}}
$$
for all $t\geq 0$ and another different positive constant $C$.
\end{Proposition}

The following proposition will be needed to study the case $p>2$. Its proof is left to the reader.

\begin{Proposition}\label{proposition:desigualdad_numerica}
For every pair of real numbers $a$ and $b$ there exists a positive finite constant $C$ such that for every $p\geq 2$
$$
(a-b)(|a|^{p-2}a - |b|^{p-2}b)\geq C(|a|^{p/2}-|b|^{p/2})^2
$$
\end{Proposition}

Thanks to these propositions we can prove the following theorem.

\begin{Theorem}
Every solution $u$ of problem~\eqref{eq:main.Cauchy.local}--\eqref{eq:main.Cauchy.nonlocal} satisfies
$$
\|u(\cdot, t)\|_{L^p(\mathbb{R}^N)}^p \leq Ct^{-\frac{N(p-1)}{2}}
$$
for every $p\in [1,\infty)$.
\end{Theorem}

\begin{Remark} {\rm
This bound coincides with the behaviour of the solutions of the local heat equation~\eqref{heat.equation}
and also with the behaviour of the solutions to the nonolocal evolution equation~\eqref{eq:nonlocal.equation}, see
\cite{Chasseigne-Chaves-Rossi-2006}.}
\end{Remark}

\begin{proof}
Inequality~\eqref{eq:fundamental.inequality} and the previous proposition provide the result when $p=2$, since $E(u)\geq \Omega_{n\ell}^J(u)$. In fact, we can just multiply by $u$ the equation and integrate to obtain
$$
\partial_t \|u\|_{L^2(\mathbb{R}^N)}^2 = - \partial_u E(u) \leq - k\int_{\mathbb{R}^N}\int_{\mathbb{R}^N} J(x-y)(u(y)-u(x))^2\, {\rm d}y{\rm d}x
$$
Using the previous proposition with $p=2$ we get
$$
\partial_t \|u\|_{L^2(\mathbb{R}^N)}^2
\leq - C
\min\Big\{\|u\|_{L^1(\mathbb{R}^N)}^{-4/N} \|u\|_{L^2(\mathbb{R}^N)}^{2(1+2/N)}, \|u\|_{L^2(\mathbb{R}^N)}^2\Big\},
$$
from where it follows that
$$
\|u(\cdot, t)\|_{L^2(\mathbb{R}^N)}^2 \leq Ct^{-\frac{N}{2}}
$$
using the conservation of mass.

The decay for $p\in(1,2)$ can be obtained through interpolation between the previous inequality and the conservation of mass property. For every $p\in(1,2)$ there exists a $\theta\in(0,1)$ such that
$$
\frac{1}{p} = \frac{1-\theta}{1} + \frac{\theta}{2}\quad\text{and}\quad \|u(\cdot,t)\|_{L^p(\mathbb{R}^N)}\leq \|u(\cdot,t)\|_{L^1(\mathbb{R}^N)}^{1-\theta}\|u(\cdot,t)\|_{L^2(\mathbb{R}^N)}^\theta.
$$
Therefore, since the mass of our solutions is constant and $\theta = 2-2/p$ we obtain that
$$
\|u(\cdot,t)\|_{L^p(\mathbb{R}^N)}\leq C\|u(\cdot,t)\|_{L^2(\mathbb{R}^N)}^{2-\frac{2}{p}}\leq C'\left( t^{\frac{-N}{4}} \right)^{2-\frac{2}{p}} = C' t^{- \frac{N(p-1)}{2p}}.
$$

Finally the case in the case $p>2$. One can check that
$$
\begin{aligned}
\partial_t \|u(\cdot,t)\|_{L^p (\mathbb{R}^N)}^p =& -\int_{\Omega_\ell}\left| \nabla (u^{p/2}(\cdot,t)) \right|^2
\\[6pt] 
&- \alpha\int_{\Omega_{n\ell}}\int_{\mathbb{R}^N}J(x-y)(u(y,t)-u(x,t))(u^{p-1}(y,t)-u^{p-1}(x,t))\,{\rm d}y{\rm d}x.
\end{aligned}$$
Using Proposition~\ref{proposition:desigualdad_numerica} and Lemma~\ref{lemma:preliminar} we arrive to
$$
\partial_t \|u(\cdot,t)\|_{L^p(\mathbb{R}^N)}^p \leq  - C\int_{\mathbb{R}^N}\int_{\mathbb{R}^N}J(x-y)(u^{p/2}(y,t)-u^{p/2}(x,t))\, {\rm d}y{\rm d}x
$$
and using~Proposition \ref{proposition:molino} with $v=u^{p/2}$ and $p=2$ (this is the $p$ from the proposition, not the $p$ of the norm we are studying) we arrive to
$$
\partial_t \|u(\cdot,t)\|_{L^p(\mathbb{R}^N)}^p \leq C\|u(\cdot,t)\|_{L^p(\mathbb{R}^N)}^{p(1+\gamma)},
$$
with $\gamma=2/N(p-1)$.
From here it is easy to finish the proof.
\end{proof}

\subsection{Rescaling the kernel}

In this part we will study how, trough a limit procedure rescaling in the kernel $J$, we can obtain the local problem.
In fact we will show that solutions to the Cauchy problem for the local heat equation~\eqref{heat.equation}
can be obtained as the limit as $\varepsilon\to 0$ of solutions $u^\varepsilon$ to the problem~\eqref{eq:main.Cauchy.local}--\eqref{eq:main.Cauchy.nonlocal} with kernel $J^\varepsilon$ given by~\eqref{eq:definition.J.epsilon} and the same initial data.

We will prove convergence of the solutions in $L^2(\mathbb{R}^N)$ for finite times with Brezis-Pazy Theorem through Mosco's convergence result and this is one of the reasons why we presented another existence of solutions result for this problem based on semi-group theory for m-accretive operators. The associated energy functional to the rescaled problem reads
$$
E^\varepsilon(u):=\int_{\Omega_\ell} \frac{|\nabla u|^2}{2} + \frac{\alpha}{2\varepsilon^{N+2}}\int_{\Omega_{n\ell}}\int_{\mathbb{R}^N} J\left(\frac{x-y}{\varepsilon}\right)(u(y)-u(x))^2\,{\rm d}y{\rm d}x
$$
if $u \in D(E^\varepsilon):=L^2(\mathbb{R}^N)\cap H^1(\Omega_\ell)$ and $E^\varepsilon(u):=\infty$ if not. Analogously, we define the limit energy functional as
$$
E'(u):=\int_{\mathbb{R}^N} \frac{|\nabla u|^2}{2}
$$
if $u \in D(E'):= H^1(\mathbb{R}^N)$ and $E'(u):=\infty$ if not.

Now, given $u_0\in L^2(\mathbb{R}^N)$, for each $\varepsilon>0$, let $u^\varepsilon$ be the solution to the evolution problem
associated with the energy $E^\varepsilon$ and initial datum $u_0$ and $u$ be the solution associated with $E'$ with the same initial condition.

\begin{Theorem}
Under the above assumptions, the functions $u^\varepsilon$ converge to solutions to \eqref{heat.equation}. For any finite $T>0$ we have that
$$
\lim\limits_{\varepsilon\to 0}\left(\sup\limits_{t\in[0,T]} \| u^\varepsilon (\cdot, t) - u(\cdot, t) \|_{L^2(\mathbb{R}^N)}\right)=0.
$$
\end{Theorem}
\begin{proof}
We will make use of the Brezis-Pazy Theorem for the sequence of m-accretive operators $B_{J^\varepsilon}$ in $L^2(\mathbb{R}^N)$ defined previously in the existence and uniqueness subsection. To apply this result we will need to show that the resolvent operators satisfy
$$
\lim\limits_{\varepsilon\to 0} (I+B_{J^\varepsilon})^{-1}\phi = (I+A)^{-1}\phi
$$
for every $\phi\in L^2(\mathbb{R}^N)$ where $A(u):=(-\Delta) u$ is the classic operator for the heat equation in this theory. If we have this then the theorem gives convergence of $u^\varepsilon$ to $u$ in $L^2(\mathbb{R}^N)$ uniformly in $[0,T]$. To prove this convergence of the resolvents we will use Mosco's result, where we only have to prove two things:
\begin{itemize}
\item[(i)] For every $u\in D(E')$ there exists a $\{u^\varepsilon\}\in D(E^\varepsilon)$ such that $u^\varepsilon\to u$ in $L^2(\mathbb{R}^N)$ and
$$
E'(u)\geq \limsup\limits_{\varepsilon\to 0} E^\varepsilon(u^\varepsilon).
$$
\item[(ii)] If $u^\varepsilon \to u$ weakly in $L^2(\mathbb{R}^N)$ then
$$
E'(u)\leq \liminf\limits_{\varepsilon\to 0} E^\varepsilon(u^\varepsilon).
$$
\end{itemize}
For more information, see~\cite{ElLibro}, Appendix 7, Theorems A.3 and A.38.

Let us start with (i). Given $u\in H^1(\mathbb{R}^N)$ we know that there exists a sequence $\{v_n\}\in C_c^\infty (\mathbb{R}^N) $ such that $v_n\to u$ in $L^2(\mathbb{R}^N)$. On the other hand, through Taylor's expansion it is not hard to see that
$$
\lim\limits_{\varepsilon\to 0} \left(\lim\limits_{\varepsilon\to 0}E^\varepsilon(v_n)\right) = E'(u).
$$
This means, by a diagonal argument, that there exists a subsequence of $\varepsilon_j$ such that
$$
\lim\limits_{\varepsilon_j\to 0}E^{\varepsilon_j}(v_j) = E'(u),
$$
showing (i).

To see (ii) from the sequence of $u^\varepsilon$ that converges weakly to $u$ we extract a subsequence $\varepsilon_n$ such that
$$
\lim\limits_{\varepsilon\to 0} E^{\varepsilon_n}(u^{\varepsilon_n}) =\liminf\limits_{\varepsilon\to 0} E^{\varepsilon_n}(u^{\varepsilon_n}).
$$
We will suppose that this inferior limit is finite, since if it is not there is nothing to prove.

Let us now take a ball of radius centered at 0, say $B^R$ and define $B^R_\ell:=B^R\cap\Omega_\ell$ and respectively $B^R_{n\ell}$ an define
$$
E_R^\varepsilon(u^\varepsilon):=\int_{B^R_\ell} \frac{|\nabla u^\varepsilon|^2}{2} + \frac{\alpha}{2\varepsilon^{N+2}}\int_{B^R_{n\ell}}\int_{B^R} J\left(\frac{x-y}{\varepsilon}\right)(u^\varepsilon(y)-u^\varepsilon(x))^2\, {\rm d}y{\rm d}x.
$$
Since the inferior limit is finite there must exist a $\varepsilon_0$ such that this quantity is bounded by a constant $M$ that only depends on $R$ for al $\varepsilon<\varepsilon_0$ and we can apply Lemma~\ref{lemma:preliminar} to obtain that there exists a positive constant $k$ not depending on $\varepsilon$ such that
$$
k\frac{\alpha}{2\varepsilon^{N+2}}\int_{B^R}\int_{B^R} J\left(\frac{x-y}{\varepsilon}\right)(u^\varepsilon(y)-u^\varepsilon(x))^2\, {\rm d}y{\rm d}x<M.
$$
Now on this domain we apply Lemma~\ref{lemma:basic_convergence_lemma_for_rescaled} to obtain a subsequence of $u^\varepsilon$, denoted by itself for simplicity, that converges to $u$ in $L^2(B^R)$ and such that
$$
\lim\limits_{\varepsilon\to 0} \left( \frac{k}{\alpha}J^{1/2}(z) \chi_{B^R}(x+\varepsilon z)\chi_{\Omega}(x) \frac{\bar{u}^\varepsilon(x+\varepsilon z) - \bar{u}^\varepsilon(x)}{\varepsilon} \right) = k\alpha J^{1/2}(z)h(x,z)
$$
weakly in $L^2_x(\mathbb{R}^N)\times L^2_z(\mathbb{R}^N)$ with $h(x,z)=z\nabla u(x)$ for all $(x,z)\in B^{R}\times \mathbb{R}^N$. Using now the lower semi-continuity of the norm for sequences that converge weakly we have that
$$
\left(\int_{\mathbb{R}^N} J(z)z^2\ dz\right) \int_{B^{R}}|\nabla u|^2
\leq \displaystyle\liminf\limits_{\varepsilon\to 0} \frac{k}{\alpha}\int_{B^{R}}\int_{B^{R}} J(x-y)\frac
{(u^\varepsilon(y,t)-u^\varepsilon(x,t))^2}{\varepsilon^2}\,{\rm d} y{\rm d}x,
$$
which means, using again Lemma~\ref{lemma:preliminar}, that
$$
\displaystyle\int_{B^{R}}|\nabla u|^2\leq \lim\limits_{\varepsilon\to 0} E_R^{\varepsilon}(u^\varepsilon)\leq \lim\limits_{\varepsilon\to 0} E^{\varepsilon}(u^\varepsilon)
$$
and we finish just by making $R$ go to $\infty$.
\end{proof}


\section{The Neumann problem} \label{sect-Neumann}
\setcounter{equation}{0}

In this section we discuss the Neumann problem~\eqref{eq:main.Neumann.local}--\eqref{eq:main.Neumann.nonlocal}.


\subsection{Existence, uniqueness and conservation of mass}

The ideas presented for the Cauchy problem can be applied \textit{mutatis mutandis} to this problem.
Therefore, using the fixed point argument, or the alternative approach by semigroup theory we obtain the following result whose
proof is left to the reader.

\begin{Theorem}\label{thm:existence_uniqueness_neumann}
Given $u_0\in L^1(\Omega)$ there exists a unique solution to problem~\eqref{eq:main.Neumann.local}--\eqref{eq:main.Neumann.nonlocal} with initial datum $u_0$. This solution conserves its mass along the evolution.
\end{Theorem}

\subsection{Comparison principle}

Also arguing as we did for the Cauchy problem we have a comparison result for the Neumann case.

\begin{Theorem}
If $u_0\geq 0$ then $u\geq 0$ in $\Omega\times\mathbb{R}_+$. Moreover, given two initial data $u_0$ and $v_0$, with $u_0\geq v_0$ then the corresponding solutions satisfy $u\geq v$ in $\Omega\times\mathbb{R}_+$.
\end{Theorem}

\subsection{Asymptotic behaviour}

In this occasion we expect the solution to converge to the average of the initial condition in every $L^p$. In fact what we are going to show is that the function
$$
v=u-|\Omega|^{-1}\int_{\Omega} u_0
$$
converges to 0 exponentially fast in $L^p$ norm.

\begin{Theorem}
The function $v$ satisfies
$$
\|v(\cdot, t)\|_{L^p(\Omega)}\leq C_1e^{-C_2 t}
$$
for every $p\in[1,\infty)$ where $C_1$ and $C_2$ are positive finite constants ($C_2$ can be taken independent of $u_0$).
\end{Theorem}

\begin{Remark}
{\rm This behaviour coincides with the behaviour of the solutions of the Heat Equation and with the 
behaviour of the solutions to the nonlocal evolution equation when the integrals are considered in the domain $\Omega$, see
\cite{Chasseigne-Chaves-Rossi-2006} (we have exponential convergence to the mean value
of the initial condition, but notice that the constants and the exponents can be different for the three cases, local, nonlocal and our mixed 
local/nonlocal problems).}
\end{Remark}

\begin{proof}
We will prove the result for $p=2$. The result for $p\in[1,\infty)$ comes from the use of H\"{o}lder's inequality and for $p>2$ from
$$
\begin{aligned}
\partial_t \|u(\cdot,t)\|_{L^p(\Omega)}^p = &-\int_{\Omega_\ell}\left| \nabla (u^{p/2})(\cdot,t) \right|^2\\[6pt] 
&- \alpha\int_{\Omega_{n\ell}}\int_{\mathbb{R}^N}J(x-y)(u(y,t)-u(x,t))(u^{p-1}(y,t)-u^{p-1}(x,t))\,{\rm d}y{\rm d}x,
\end{aligned}
$$
the Proposition~\ref{proposition:desigualdad_numerica}, Lemma~\ref{lemma:preliminar} and the fact that
$$
\nabla u\cdot \nabla(|u|^{p-1}u) = |\nabla (|u|^{p/2})|^2
$$
so we can rename $u^{p/2}$ as another function $w$ and apply the case $p=2$ to obtain the cases $p>2$. This case is again left for the reader
and is somehow similar to the analogous case for the Cauchy problem.

So for $p=2$ we compute, after some calculations
$$
\partial_t\int_{\Omega}v^2(x,t) \, {\rm d}x= -2\int_{\Omega_\ell}|\nabla v|^2 - 2C\int_{\Omega_{n\ell}}\int_\Omega J(x-y)(v(y,t)-v(x,t))^2\, {\rm d}y{\rm d}x.
$$
Using Lemma~\ref{lemma:preliminar} and a result from \cite{ElLibro} that shows that for every function $v$ with zero average we have that
$$
\int_{\Omega}\int_\Omega J(x-y)(v(y,t)-v(x,t))^2\, {\rm d}y{\rm d}x \geq \beta \int_{\Omega}v^2 (x,t) \, {\rm d}x
$$
for some positive $\beta$ we obtain
$$
\partial_t\int_{\Omega}v^2 (x,t) \, {\rm d}x \leq -k\int_{\Omega}v^2 (x,t) \, {\rm d}x
$$
for another constant $k$. From this point the result follows trivially.
\end{proof}

\begin{Remark} {\rm
One can define what is the analogous to the first non-zero eigenvalue for this problem as
$$
\beta_1 (\Omega) = \inf_{\int_\Omega u = 0} \frac{ \displaystyle \int_{\Omega_\ell} |\nabla u|^2 + \alpha
\int_{\Omega_{n\ell}}\int_{\Omega} J(x-y)(u(y)-u(x))^2\, {\rm d}y{\rm d}x}{ \displaystyle  \int_\Omega u^2 }
$$
and show that $\beta_1 (\Omega)$ is positive (otherwise the exponential decay does not hold
with a constant $C_2$ independent of $u_0$). The existence of an eigenfunction (a function that achieves the
infimum) is not straightforward. We leave this fact open.
}
\end{Remark}

\subsection{Rescaling the kernel}

As before, we want to study the convergence  of the solutions $u^\varepsilon$ of the Neumann  problem~\eqref{eq:main.Neumann.local}--\eqref{eq:main.Neumann.nonlocal} with rescaled kernel  $J^\varepsilon$ given by~\eqref{eq:definition.J.epsilon} to the solution of the Neumann problem for the local heat equation with the same initial datum.

\begin{Theorem}\label{thm:rescaled_Neumann}
For any finite $T>0$ we have that
$$
\lim\limits_{\varepsilon\to 0}\left(\sup\limits_{t\in[0,T]} \| u^\varepsilon (\cdot, t) - u(\cdot, t) \|_{L^2(\mathbb{R}^N)}\right)=0
$$
where this $u$ is the solution of the Neumann problem for the Heat Equation in $\Omega$ with the same initial data $u_0$.
\end{Theorem}

The proof of this theorem is analogous to the one we did for the Cauchy problem (see also~\cite{BLGneu}). Again we use the already mentioned Brezis-Pazy Theorem through convergence of the resolvents. Notice that this approach uses the linear semi-group theory in $L^2(\Omega)$ mentioned in the Cauchy section (that also works just fine in this case).


\section{The Dirichlet problem} \label{sect-Dirichlet}
\setcounter{equation}{0}

In this section we devote our attention to the Dirichlet problem~\eqref{eq:main.Dirichlet.local}--\eqref{eq:main.Dirichlet.boundary}.


\subsection{Existence and uniqueness}

Again, we have the following result whose proof can be obtained as in the previous cases (again
we have two proofs, one using a fixed point argument and another one using semigroup theory).

\begin{Theorem}\label{thm:existence_uniqueness_dirichlet}
Given $u_0\in L^1(\Omega)$, there exists a unique solution to problem~\eqref{eq:main.Dirichlet.local}--\eqref{eq:main.Dirichlet.boundary} which has $u_0$ as initial datum.
\end{Theorem}

\begin{Remark} {\rm For this problem there is a loss of mass trough the boundary. In fact, assume that $u_0$ is nonnegative
(and hence $u(x,t)$ is nonnegative for every $t>0$). Then integrating in $\Omega$ we get
$$
\begin{aligned}
\partial_t \int_\Omega u(x,t) \, {\rm d}x =& \int_{\partial \Omega \cap \partial \Omega_\ell}
\partial_\eta u(x,t) \, {\rm d}\sigma  +  \alpha \int_{\Omega_\ell} \int_{\mathbb{R}^N\setminus \Omega_\ell} J(x-y)(u(y,t)-u(x,t))\, {\rm d}y{\rm d}x
\\[6pt]&+ \alpha\int_{\Omega_{n\ell}} \int_{\mathbb{R}^N\setminus \Omega_{n\ell}} J(x-y)(u(y,t)-u(x,t))\,{\rm d}y{\rm d}x
\\[6pt]
=& \int_{\partial \Omega \cap \partial \Omega_\ell}
\frac{\partial u}{\partial \eta}(x,t) \, d\sigma(x) -  \alpha \int_{\Omega_\ell} \int_{\mathbb{R}^N\setminus \Omega} J(x-y)u(x,t) \, {\rm d}y{\rm d}x
\\[6pt]
& - \alpha  \int_{\Omega_{n\ell}} \int_{\mathbb{R}^N\setminus \Omega} J(x-y)u(x,t) \, {\rm d}y{\rm d}x <0.
\end{aligned}
$$
}
\end{Remark}

\subsection{Comparison principle}

As for the previous cases we have a comparison result.

\begin{Theorem}
If $u_0\geq 0$ then $u\geq 0$ in $\Omega\times\mathbb{R}_+$. Moreover, given two initial data $u_0$ and $v_0$ with $u_0\geq v_0$,  then the corresponding solutions satisfy
$u\ge v$ in $\Omega\times\mathbb{R}_+$.
\end{Theorem}

\subsection{Asymptotic decay}

The result here is analogous to the one in corresponding section for the Neumann problem. The only extra tool needed is a result that was proved in \cite{ElLibro} that shows that there exists a constant $\beta$ such that for every function $u$ that satisfies $u(x,t)\equiv 0$ for every $x\in\mathbb{R}^N\setminus \Omega$ we have that
$$
\int_{\mathbb{R}^N}\int_\Omega J(x-y)(u(y,t)-u(x,t))^2\,{\rm d}y{\rm d}x \geq \beta \int_{\Omega}u^2 (x)\, {\rm d}x ,
$$
similarly to the previous section for functions with zero average. At this point the proof for the following theorem is straightforward.

\begin{Theorem}
The solution $u$ of the Dirichlet problem problem~\eqref{eq:main.Dirichlet.local}--\eqref{eq:main.Dirichlet.boundary} with initial datum $u_0$  satisfies
$$
\|u(\cdot, t)\|_{L^p(\Omega)}\leq C_1e^{-C_2 t}
$$
for every $p\in[1,\infty)$ where $C_1$ and $C_2$ are positive finite constants ($C_2$ can be chosen independent of $u_0$).
\end{Theorem}

\begin{Remark} {\rm
This coincides with the behaviour of the solutions to the Heat Equation and with the 
behaviour of the solutions to the corresponding nonlocal evolution equation with zero exterior condition, see
\cite{Chasseigne-Chaves-Rossi-2006} (we have exponential decay, but notice that again here 
the constants and the exponents can be different for the three cases).}
\end{Remark}

\begin{Remark} {\rm Again for this case we have an associated eigenvalue problem. Let us consider
$$
\lambda_1 (\Omega) = \inf_{u|_{\mathbb{R}^N \setminus \Omega} \equiv 0} \frac{ \displaystyle \int_{\Omega_\ell} |\nabla u|^2 + \alpha
\int_{\Omega_{n\ell}}\int_{\mathbb{R}^N} J(x-y)(u(y)-u(x))^2\, {\rm d}y{\rm dx}}{ \displaystyle  \int_\Omega u^2}.
$$
One can prove that $\lambda_1 (\Omega)$ is positive using our control of the nonlocal energy, \eqref{eq:fundamental.inequality}
and the results in \cite{ElLibro}. Again in this case the existence of an eigenfunction (a function that achieves the
infimum) is left open.
}
\end{Remark}

\subsection{Rescaling the kernel}

In this part we will study how we can obtain the solution to the Dirichlet problem with zero boundary datum for the local heat equation as the limit as $\varepsilon\to 0$ of solutions $u^\varepsilon$  of the Dirichlet problem~\eqref{eq:main.Dirichlet.local}--\eqref{eq:main.Dirichlet.boundary} with rescaled kernel  $J^\varepsilon$ as in~\eqref{eq:definition.J.epsilon} and the same initial datum.
\begin{Theorem}
For any finite $T>0$ we have that
$$
\lim\limits_{\varepsilon\to 0}\left(\sup\limits_{t\in[0,T]} \| u^\varepsilon (\cdot, t) - u(\cdot, t) \|_{L^2(\mathbb{R}^N)}\right)=0
$$
where this $u$ is the solution of the Dirichlet problem for the Heat Equation in $\Omega$ with the same initial data $u_0$ and zero boundary data.
\end{Theorem}

The proof of this theorem is analogous to the previous ones (we also refer to~\cite{BLGdir} here), see the comments 
about Theorem~\ref{thm:rescaled_Neumann} in
the previous section.

\section{Comments on possible extensions}
\setcounter{equation}{0}

In this final section we briefly comment on possible extensions of our results.

\begin{itemize}
\item Our results could be extended to cover singular kernels including, for example, fractional Laplacians.
In this case the associated energy for the Cauchy problem looks like
$$
E(u):=\int_{\Omega_\ell} \frac{|\nabla u|^2}{2}+
\frac{C}{2}\int_{\mathbb{R}^N_-}\int_{\mathbb{R}^N} \frac{(u(y)-u(x))^2}{|x-y|^{N+2s}}
\, {\rm d}y{\rm d}x.
$$
The abstract semigroup theory seems the right way to obtain existence and uniqueness of a solution.
One interesting problem is to couple two different fractional Laplacians and look for the asymptotic behaviour
of the solutions to the corresponding Cauchy problem. We will tackle this kind of extension
in a future paper.

\item One can look for moving interfaces, making that $\Gamma$ depends on $t$. To show existence and uniqueness
of solutions for a problem like this seems a challenging problem. In this framework one is tempted to consider free boundary problems
in which we have an unknown interface that evolves with time and we impose that solutions have conservation of the total mass plus some continuity  across the free boundary.

\item Finally, we mention that an interesting problem is to look at nonlinear diffusion equations (coupling, for example, a
local $p-$Laplacian with a nonlocal $q-$Laplacian, see \cite{ElLibro} for a definition of the last operator). A possible
energy for this problem is
$$
E(u):=\int_{\Omega_\ell} \frac{|\nabla u|^p}{p} +
\frac{C}{q}\int_{\Omega_{n\ell}}\int_{\mathbb{R}^N} \frac{|u(y)-u(x)|^q}{|x-y|^{N+qs}}
\,{\rm d}y{\rm d}x.
$$
This problem involves new difficulties, especially when one looks
for scaled versions of the kernel and tries to see whether there is a limit.

\end{itemize}

\medskip

\noindent{\large \textbf{Acknowledgments}}

\noindent The first two authors were partially supported by the Spanish project {MTM2017-87596-P}, and the third one by CONICET grant PIP GI No 11220150100036CO
(Argentina), by  UBACyT grant 20020160100155BA (Argentina) and by the Spanish project MTM2015-70227-P.

This project has received funding from the European Union’s Horizon 2020 research and innovation programme under the Marie Skłodowska-Curie grant agreement No.\,777822.

Part of this work was done during visits of JDR to Madrid and of AG and FQ to Buenos Aires. The authors want to thank these institutions for the nice and stimulating working atmosphere.



\

\noindent\textbf{Addresses:}

\noindent\textsc{A. G\'arriz: } Departamento de Matem\'{a}ticas, Universidad
Aut\'{o}noma de Madrid, 28049 Madrid, Spain. (e-mail: alejandro.garriz@estudiante.uam.es).

\noindent\textsc{F. Quir\'os: } Departamento de Matem\'{a}ticas, Universidad
Aut\'{o}noma de Madrid, 28049 Madrid, Spain. (e-mail: fernando.quiros@uam.es).

\noindent\textsc{J. D. Rossi: } Departamento de Matem\'{a}ticas, FCEyN, Universidad
de Buenos Aires, Ciudad Universitaria, Pab. 1 (1428) Buenos Aires, Argentina. (e-mail: jrossi@dm.uba.ar).

\end{document}